\documentclass{amsart}
\usepackage{amsmath, amsthm, amsfonts, hyperref, graphicx, ifpdf}
\usepackage{esint}
\usepackage[dvipsnames,usenames]{color}
\hypersetup{
   unicode=false,          
   pdftoolbar=true,        
   pdfmenubar=true,        
   pdffitwindow=false,     
   pdfstartview={},    
   pdftitle={},    
   pdfauthor={},     
   pdfsubject={},   
   pdfcreator={},   
   pdfproducer={}, 
   pdfkeywords={}, 
   pdfnewwindow=true,      
   colorlinks=true,       
   linkcolor=blue,          
   citecolor=blue,        
   filecolor=blue,      
   urlcolor=blue          
}

\newcommand{\R}{{\mathbb R}}

\newtheorem{theorem}{Theorem}[section]
\newtheorem{lemma}[theorem]{Lemma}
\newtheorem{pro}[theorem]{Proposition}
\newtheorem{cor}[theorem]{Corollary}

\theoremstyle{definition}
\newtheorem{Def}[theorem]{Definition}

\theoremstyle{remark}
\newtheorem{remark}[theorem]{Remark}

\numberwithin{equation}{section}

\newcommand{\cmc}{constant mean curvature}
\newcommand{\es}{Euclidean space}
\newcommand{\ee}{evolution equation}
\newcommand{\evs}{evolving surface}
\newcommand{\hmc}{h-mean convex}
\newcommand{\hs}{hyperbolic space}
\newcommand{\mcv}{mean convex}
\newcommand{\mc}{mean curvature}
\newcommand{\mcf}{mean curvature flow}
\newcommand{\maxp}{maximum principle}
\newcommand{\pc}{principal curvature}

\newcommand{\sapmcf}{surface area preserving mean curvature flow}
\newcommand{\sff}{second fundamental form}

\newcommand{\ti}{time interval}

\newcommand{\vpmcf}{volume preserving mean curvature flow}
\newcommand{\wrt}{with respect to}

\newcommand{\be}{\begin{equation}}
\newcommand{\ene}{\end{equation}}
\newcommand{\br}{\begin{remark}}
\newcommand{\er}{\end{remark}}
\newcommand{\bl}{\begin{lemma}}
\newcommand{\el}{\end{lemma}}
\newcommand{\bcor}{\begin{cor}}
\newcommand{\ecor}{\end{cor}}
\newcommand{\bpro}{\begin{pro}}
\newcommand{\epro}{\end{pro}}
\newcommand{\ben}{\begin{enumerate}}
\newcommand{\een}{\end{enumerate}}
\newcommand{\bp}{\begin{proof}}
\newcommand{\ep}{\end{proof}}
\newcommand{\bpo}{\begin{pro}}
\newcommand{\epo}{\end{pro}}
\newcommand{\beq}{\begin{equation*}}
\newcommand{\eeq}{\end{equation*}}
\newcommand{\bear}{\begin{eqnarray}}
\newcommand{\eear}{\end{eqnarray}}
\newcommand{\beqar}{\begin{eqnarray*}}
\newcommand{\eeqar}{\end{eqnarray*}}
\newcommand{\bt}{\begin{theorem}}
\newcommand{\et}{\end{theorem}}

\DeclareMathAlphabet{\mathcal}{OMS}{cmsy}{m}{n}

\numberwithin{equation}{section}

\allowdisplaybreaks

\def\XXint#1#2#3{{\setbox0=\hbox{$#1{#2#3}{\int}$}
    \vcenter{\hbox{$#2#3$}}\kern-.5\wd0}}

\makeatletter
\def\@citestyle{\m@th\upshape\mdseries}
\def\citeform#1{{\bfseries#1}}
\def\@cite#1#2{{%
  \@citestyle[\citeform{#1}\if@tempswa, #2\fi]}}
\@ifundefined{cite }{%
  \expandafter\leqt\csname cite \endcsname\cite
  \edef\cite{\@nx\protect\@xp\@nx\csname cite \endcsname}%
}{}
\makeatother

\newcommand{\ppl}[2]{\frac{\partial{#1}}{\partial{#2}}}

\renewcommand{\H}{\mathbb{H}}

\begin{document}
\newcommand{\osc}{{\text{osc}}}
\newcommand{\Vol}{{\text{Vol}}}
\newcommand{\V}{{\text{V}}}
\newcommand{\nn}{{\text{n}}}
\newcommand{\cM}{{\cal{M}}}
\newcommand{\Ric}{{\text{Ric}}}
\newcommand{\RE}{{\text{Re }}}
\newcommand{\LL}{{\cal{L}}}
\newcommand{\diam}{{\text {diam}}}
\newcommand{\dist}{{\text {dist}}}
\newcommand{\Area}{{\text {Area}}}
\newcommand{\Length}{{\text {Length}}}
\newcommand{\Energy}{{\text {Energy}}}
\newcommand{\SSS}{{\bold S}}
\newcommand{\K}{{\text{K}}}
\newcommand{\Hess}{{\text {Hess}}}
\def\RR{{\bold  R}}
\def\SS{{\bold  S}}
\def\TT{{\bold  T}}
\def\CC{{\bold C }}
\newcommand{\dv}{{\text {div}}}
\newcommand{\kg}{{\text {k}}}
\newcommand{\expp}{{\text {exp}}}
\newcommand{\e}{{\text {e}}}
\newcommand{\eqr}[1]{(\ref{#1})}
\newcommand{\ec}{\varepsilon_c}
\newcommand{\rc}{\rho_c}
\newcommand{\sztp}{\Sigma^{0,2\pi}}

\def\p{\partial}

\parskip1ex






\title[stability of VPMCF in $\H^{n+1}$]{Stability of the Volume Preserving Mean Curvature Flow in Hyperbolic Space}

\author{Zheng Huang}
\address[Z. ~Huang]{Department of Mathematics, The City University of New York, Staten Island, NY 10314, USA}
\address{The Graduate Center, The City University of New York, 365 Fifth Ave., New York, NY 10016, USA}
\email{zheng.huang@csi.cuny.edu}

\author{Longzhi Lin}
\address[L.~Lin]{Mathematics Department\\University of California, Santa Cruz\\1156 High Street\\
Santa Cruz, CA 95064\\USA}
\email{lzlin@ucsc.edu}

\author{Zhou Zhang}
\address[Z. ~Zhang]{School of Mathematics and Statistics, The University of Sydney, NSW 2006, Australia}
\email{zhangou@maths.usyd.edu.au}



\begin{abstract}
We consider the dynamic property of the {\vpmcf}. This flow was introduced by Huisken in \cite{Hui87} who also proved it converges to a round sphere of the same enclosed volume if the initial hypersurface is strictly convex in Euclidean space. We study the stability of this flow in hyperbolic space. In particular, we prove that if the initial hypersurface is hyperbolically {\mcv} and close to an umbilical sphere in the $L^2$-sense, then the flow exists for all time and converges exponentially to an umbilical sphere.

\end{abstract}

\maketitle

\section{Introduction}
\subsection{Background and Main Theorem}
Let $M^n$ be a smooth, embedded, closed (compact, no boundary) $n$-dimensional manifold in {\hs} $\H^{n+1}$ ($n \ge 2$), and we evolve it by the {\vpmcf} (VPMCF),
\be
\label{eq-vp}
\frac{\partial F}{\partial t} = \left(h - H\right)\nu, \qquad F(\cdot,0) = F_0(\cdot)\,,
\ene
where $F_0: M^n \to \H^{n+1}$ is the initial embedding, $H = H(x,t)$ is the {\mc} and $\nu = \nu (x,t)$ is the outward unit normal vector of the {\evs} 
$M_t = F(\cdot, t)$ at point $(x,t)$ (for simplicity, we write $(x,t)\in M_t$). The function $h$ is the average of the {\mc} on $M_t$, given by 
\be\label{def-h}
h = h (t) = \fint_{M_t}H\, d\mu =\frac{\int_{M_t} H \, d\mu}{\int_{M_t} \, d\mu} \,,
\ene
where $d\mu = d\mu_t$ denotes the surface area element of the {\evs} $M_t$ {\wrt} the induced metric $g(t)$. 

In this paper we use the convention that the {\mc} is the sum of all {\pc}s. Clearly we have $H\not\equiv 0$ on $M_0$ since there is no closed minimal 
hypersurface in {\hs}. The presence of the global term $h$ in the VPMCF equation \eqref{eq-vp} forces the flow to behave quite differently from the usual {\mcf} (MCF).

Hypersurfaces of {\cmc} are critical points of the area functional under the constraint of fixed enclosed volume. These
hypersurfaces are also static state for the VPMCF equation \eqref{eq-vp}. A remarkable theorem of Huisken-Yau (\cite{HY96})
on the existence of a foliation of spheres outside of some large compact set in asymptotic flat manifolds was achieved by
studying a parameter family of VPMCFs. This flow, and the {\sapmcf} studied in \cite{McC03, HL15}, are special cases of so-called 
mixed volume preserving {\mcf}. They are closely related to convex geometry and classical inequalities, see for instance 
\cite{ES98, McC04, ACW21, WX14}, also see \cite{CRM12, EM12} for other geometric settings.

We denote $A = \{a_{ij}\}$ the {\sff} of $M_t$ and ${\text\AA}= A - \frac{H}{n}g$ its traceless part. Then we have
$|{\text\AA}|^2 = |A|^2 - \frac{1}{n}H^2$. This quantity measures the roundness of a (closed, immersed) hypersurface 
$\Sigma$ in $\mathbb{H}^{n+1}$: if ${\text\AA}\equiv 0$, i.e., $\Sigma$ is umbilic at every point, then by a classical 
Codazzi's theorem in differential geometry, it is a geodesic sphere, see e.g. \cite[Theorem 29]{Spi79}. We also remark 
that, in $\R^3$, De Lellis and M{\"u}ller \cite{DLM05} generalized Codazzi's theorem by showing a version of the 
following quantitative rigidity
\beq
\inf_{\lambda\in \mathbf{R}} \| A - \lambda \,\text{Id}\|_{L^2(\Sigma)} \leq C \|{\text\AA}\|_{L^2(\Sigma)}\,,
\eeq
for some universal constant $C$. Such quantitative rigidity is not available for {\hs} to our acknowledge.

Strict convexity (i.e., all {\pc}s are positive) plays a fundamental role in classical works of several types of MCFs, especially 
in {\es}. Huisken (\cite{Hui84}) proved that an initial smooth closed and strictly convex hypersurface will stay convex and 
flow into a round point along the MCF in Euclidean space. He (\cite{Hui87}) also showed, in the case of the VPMCF, the 
flow of an initial smooth closed and strictly convex hypersurface will exist for all time and flow into a round sphere in Euclidean space. The 
parallel result for the {\sapmcf} is also true, showed by McCoy (\cite{McC03}). Though natural in Euclidean geometry, 
this notion of convexity is not the most natural in {\hs}. The presence of horospheres in {\hs} poses strong restrictions 
on the geometry of hypersurfaces (via Hopf's {\maxp}): for instance any closed {\cmc} hypersurface has {\mc} greater 
than $n$ in $\H^{n+1}$. 

\begin{Def}\label{def}
We call a hypersurface of an $(n+1)$-dimensional hyperbolic manifold (strictly) {\it h-convex} if every {\pc} of the hypersurface at every point is greater than $1$, and call it (strictly) {\it h-{\mcv}} or hyperbolically {\mcv} if its {\mc} at every point is greater than $n$. 
\end{Def}

The ``h-convexity" was introduced in (\cite{CRM07}), where the authors proved that h-convexity is preserved along the VPMCF
in {\hs}. Moreover, under the assumption of closed initial hypersurface being h-convex, they showed that the {\vpmcf} exists for all time 
and converges to an umbilical sphere. The ``h-mean convexity", or the notion of being hyperbolically {\mcv}, is much weaker than h-convexity, 
and it is not known to be preserved along the VPMCF. But it turns out this condition plays a very important role in proving the dynamic stability of the VPMCF.

Unlike the regular MCF, the VPMCF  \eqref{eq-vp} has a global forcing term in the equation which greatly complicates the analysis of the flow. 
How the singularities of the flow may form remains elusive at the current stage of study, even in Euclidean space. Moreover in our {\hs} setting, the 
negative curvature of the ambient space presents significant challenges in analyzing the {\ee}s involved in the study. As a first step to 
understand the long term behavior of the flow, in this paper, we study the dynamical property of the VPMCF \eqref{eq-vp} in {\hs} 
in the situation that the initial hypersurface is not necessarily h-convex, yet close to an umbilical sphere in the $L^2$-sense. More precisely, 
we show the stability of the flow with initial {\it h-mean convex} hypersurface (namely the initial {\mc} at every point is greater than $n$) and small 
$L^2$-norm of the traceless part of the {\sff}. Our main theorem is the following:
\bt\label{main}
Let $M_t^n \subset \H^{n+1}, n \geq 2,$ be a smooth closed solution to the VPMCF \eqref{eq-vp} for
$t\in [0,T)$ with $T\leq \infty$. Assume that $M_0$ is h-{\mcv} and
\be\label{Condition-1}
\max\left\{|M_0|^2\,, \, \max_{M_0}|A|^2\,, \,\int_{M_0}|\nabla^m A|^2\,d\mu \right\}\,\leq \, \Lambda_0^2, 
\ene
for some $\Lambda_0 \gg 1$ and all $m \in [1, n+3]$, where $|M_t|$ is the 
$n$-dimensional surface area of $M_t$ with the induced metric. Then there exists some 
$\epsilon_0 = \epsilon_0(n,\Lambda_0)>0$ such that if
\be\label{Condition-2}
\int_{M_0}|{\text\AA}|^2 \,d\mu < \epsilon_0\,,
\ene
then $T = \infty$ and the flow converges exponentially to an umbilical sphere which encloses the same volume as $M_0$.
\et
\br
It is very important that $\epsilon_0$ in Theorem \ref{main} does not depend on the lower bound of $H - n>0$ on $M_0$ from $h$-mean convexity.
\er

\subsection{Outline of the proof:}
We would like to stress that there are several serious complications in order to investigate the dynamic stability for VPMCF: with a forcing term $h$, the flow is global in nature, therefore it's difficult to localize the analysis and it is essential to keep track of $h$ along the flow; we are working in the {\hs} where $h$-{\mcv}ity is likely not preserved along the flow in general and the negative curvature of hyperbolic space makes the analysis of the flow substantially more involved. To overcome these difficulties, we use iteration techniques in combination of 
several new tools to prove the main theorem.

We will organize the iteration argument in four steps: step 1, based on the initial bounds, we derive bounds on some short time interval for several 
geometric quantities (Lemmas \ref{key} and \ref{key2}) such as $H$, $\nabla H$, $|{\text\AA}|^2$, etc. As a consequence, we show that
the {\hmc}ity is preserved on some definite time interval provided the initial hypersurface is close enough to an umbilical sphere in the $L^2$ sense; in step 2, we prove exponential decay for these quantities on the time 
interval obtained in previous step (Theorem ~\ref{exp1}), which allows us to obtain uniform bounds for these quantities on the 
interval; in step 3, we repeat above arguments to extend the time interval (Theorem ~\ref{ext}), and finally in step 4, we prove 
the amount of extension for time only depends on the initial conditions. Main theorem then follows.



\subsection{Acknowledgements} 
The research of Z. Huang is partially supported by a PSC-CUNY grant. The research of Z. Zhang is partially supported by ARC Future Fellowship FT150100341. 

\section{Technical Preparations}

In this section, we collect basic {\ee}s for key quantities, and derive some preliminary estimates that will be used in the proof. 

\subsection{Evolution equations}

Let us first fix some notations of the following geometric quantities that will be used in this study:
\ben
\item
the induced metric of the evolving hypersurface $M_t$: $\{g_{ij}(t)\}$;
\item
the {\sff} of $M_t$: $A(\cdot,t) = \{a_{ij}(\cdot,t)\}$, and its square norm is given by
\beq
|A(\cdot,t)|^2 = g^{ij}g^{kl}a_{ik}a_{jl};
\eeq
\item
the {\mc} of $M_t$ {\wrt} the outward unit normal vector: $H(\cdot,t) = g^{ij}a_{ij}$;
\item
the traceless part of the {\sff}: ${\text\AA} = A - \frac{H}{n}g$;
\item
the area element of the evolving hypersurface $M_t$: $d\mu_t = \sqrt{det(g_{ij})}$.
\een
The {\ee}s for these quantities are as follows:
\bl (\cite{Hui87, HY96})
The metric of $M_t$ satisfies the {\ee}
\be\label{eq-g}
\frac{\partial}{\partial t} g_{ij} = 2(h-H)a_{ij}\,.
\ene
Therefore,
\be
\frac{\partial}{\partial t} g^{ij} = -2(h-H)a^{ij}
\ene
and the area element satisfies:
\be\label{mu}
\frac{\partial}{\partial t} (d\mu_t) = H(h-H)d\mu_t.
\ene
Moreover, the outward unit normal vector $\nu$ to $M_t$ satisfies
\be
\frac{\partial \nu}{\partial t}  = \nabla H\,,\, 
\ene
where $\frac{\partial \nu}{\partial t}$ is a conventional way of writing down $\bar\nabla_{\frac{\p}{\p t}}\nu$.
\el
By \eqref{mu}, we have the following geometrical properties of the VPMCF:
\bcor (\cite{Hui87})\label{HuiLemma2}
\ben
\item
The $(n+1)$-dimensional volume $V_t$ of the region enclosed by $M_t$ remains unchanged along the flow, i.e.,
\beq
\frac{d}{dt}V_t = \int_{M_t}(h-H)\,d\mu = 0\,.
\eeq
\item
The $n$-dimensional surface area $|M_t|$ of $M_t$ is non-increasing along the flow, i.e.,
\beq
\frac{d}{dt}|M_t| = \frac{d}{dt}\int_{M_t} \,d\mu = \int_{M_t} H(h-H)\,d\mu = -\int_{M_t} (h-H)^2\,d\mu \le 0\,.
\eeq
\een
\ecor

Following Huisken's calculations for the MCF in general Riemannian manifolds (\cite{Hui86}), we have the
following  {\ee}s for key quantities in our setting. See also \cite{HY96} for the case $n=2$ in the setting of
asymptotic flat manifolds and \cite{CRM07} for equivalent formulas in {\hs} setting. 

\bt \label{Evol-A-H}
We have the {\ee}s for $H$ and $|A|^2$: 

\begin{itemize}
\item[(i)]
\be\label{eq-H}
\frac{\partial}{\partial t} H = \Delta H + (H-h)(|A|^2-n); 
\ene
\item[(ii)]
\be\label{eq-A2}
\frac{\partial}{\partial t} |A|^2 = \Delta |A|^2 - 2|\nabla A|^2 + 2|A|^2(|A|^2+n) - 2h\,\text{tr}\left(A^3\right) + 2H(h-2H).
\ene
\end{itemize}
where $\text{tr}\left(A^3\right) = g^{ij}g^{kl}g^{mn}a_{ik}a_{lm}a_{nj}$.
\et
We include a short proof for readers' convenience.
\bp
Let $\bar{g} = \{\bar{g}_{\alpha\beta}\}$ be the metric on $\H^{n+1}$, $\bar{\nabla}$ and $\bar{R}_{\alpha\beta\gamma\delta}$ be covariant 
derivative and Riemannian curvature tensor {\wrt} $\bar{g}$. The equation \eqref{eq-H} is clear since $\overline{\rm Ric}(\nu, \nu)=-n$ in $\H^{n+1}$. 
For \eqref{eq-A2}, we first follow \cite{Hui86, CRM07} to find that the {\sff} $\{a_{ij}\}$ of $M_t$ satisfies
the following {\ee}: 
\bear\label{eq-a}
\frac{\partial}{\partial t}a_{ij} &=& \Delta a_{ij} +(h- 2H)a_{i\ell}a_{j\ell} +|A|^2a_{ij} + a_{ij}\bar{R}_{0\ell0\ell}
-h\bar{R}_{0i0j} \notag\\
&\ & -a_{j\ell}\bar{R}_{\ell mim}-a_{i\ell}\bar{R}_{\ell mjm}+2a_{\ell m}\bar{R}_{\ell imj}
-\bar{\nabla}_j \bar{R}_{0\ell i\ell} -\bar{\nabla}_{\ell} \bar{R}_{0ij\ell}. 
\eear
The last two terms which involve the covariant derivatives of the curvature tensor drop out as we are in a constant 
curved space. Furthermore, since $\H^{n+1}$ has constant sectional curvature $-1$, the Riemannian curvature tensor 
is given by:
\be\label{rie-tensor}
\bar{R}_{\alpha\beta\gamma\delta} = (-1)\cdot (\bar{g}_{\alpha\gamma}\bar{g}_{\beta\delta}-
\bar{g}_{\alpha\delta}\bar{g}_{\gamma\beta}).
\ene
Now \eqref{eq-A2} follows from contraction and \eqref{eq-g}.
\ep

The covariant derivatives for $A$ satisfy the following.

\bcor\label{hdA-1}
We have the {\ee} for $|\nabla^m A|^2$ with $m\geq 1$:
\begin{align}\label{hdA-2}
\frac{\partial}{\partial t} |\nabla^m A|^2 = & \Delta |\nabla^m A|^2 - 2|\nabla^{m+1} A|^2 + \nabla^m A\ast \nabla^m A
\notag\\
& +\sum_{i+j+k=m} \nabla^i A\ast \nabla^j A \ast \nabla^k A \ast \nabla^m A + \sum_{r+s=m} 
\nabla^r A\ast \nabla^s A \ast \nabla^m A  \,,
\end{align}
where $S\ast \Omega$ denotes any linear combination of tensors formed by contraction on $S$ and $\Omega$ 
by the metric $g$. Here, in addition to constants, $h = h(t)$ (having only time variable) may be involved in the coefficients of the contraction.
\ecor

\bp

We have the following evolution of the second fundamental form from the proof of Theorem \ref{Evol-A-H}: 
\be
\label{eq-A}
\frac{\partial}{\partial t}A = \Delta A + A\ast A\ast A+A\ast A+ \ast A. \notag
\ene
Meanwhile, the time derivative of the Christoffel symbols $\Gamma_{jk}^i$ is equal to 
\begin{align}\label{Chris}
\frac{\partial}{\partial t} \Gamma_{jk}^i &= \frac{1}{2}g^{il} \left\{\nabla_j\left(\frac{\partial}{\partial t} g_{kl}\right) +
\nabla_k\left(\frac{\partial}{\partial t} g_{jl}\right)-\nabla_l\left(\frac{\partial}{\partial t} g_{jk}\right)\right\} \notag \\
&= g^{il} \left\{\nabla_j\left( (h- H)a_{kl}\right) + \nabla_k\left((h- H)a_{jl}\right)-\nabla_l\left((h- H)a_{jk}\right)\right\} \\
&= *\nabla A + A \ast \nabla A\,, \notag
\end{align}
where $\ast T$ denotes contraction of $T$ by the metric $g$. Note that we have used the {\ee} \eqref{eq-g} for the metric. 

Now we proceed as in \cite[\S 13]{Ham82} (see also \cite[\S 7]{Hui84}) to obtain the equation \eqref{hdA-2}. In particular, 
using \eqref{Chris}, if $S$ and $\Omega$ are tensors satisfying the {\ee} $\frac{\partial}{\partial t} S = \Delta S + \Omega$, 
then the covariant derivative $\nabla S$, which involves the Christoffel symbols, satisfies an equation of the following form:
\be\label{Covariant-1}
\frac{\partial}{\partial t} \nabla S = \Delta (\nabla S) + S*A \ast \nabla A + S*\nabla A + A\ast A\ast \nabla S + \nabla \Omega.
\ene
Therefore by \eqref{eq-a}, we find 
\beq
\frac{\partial}{\partial t} \nabla A = \Delta \nabla A  +
\sum_{i+j+k=1} \nabla^i A\ast \nabla^j A \ast \nabla^k A + \sum_{r+s=1} \nabla^r A\ast \nabla^s A + \ast\nabla A \,.
\eeq
Then by induction we have for $m\geq 1$,
\be \label{eq-a3}
\frac{\partial}{\partial t} \nabla^m A = \Delta \nabla^m A  +
\sum_{i+j+k=m} \nabla^i A\ast \nabla^j A \ast \nabla^k A + \sum_{r+s=m} \nabla^r A\ast \nabla^s A + \ast\nabla^m A \,.
\ene
Then the equation \eqref{hdA-2} follows from the following identity essentially from \eqref{eq-g}\beq
\frac{\partial}{\partial t} |\nabla^m A|^2 = 2 \left\langle \nabla^m A, \frac{\partial}{\partial t} \nabla^m A\right\rangle + 
A \ast\nabla^m A\ast\nabla^m A + A \ast A \ast\nabla^m A\ast\nabla^m A
\eeq
and the standard identity
\beq
\Delta |\nabla^m A|^2 = 2\left\langle \nabla^m A, \Delta \nabla^m A\right\rangle + 2 |\nabla^{m+1} A|^2.
\eeq
\ep

We also have the time derivative for the average of {\mc} $h(t)$.
\bl\label{h-der}
\be\label{derivativeh}
h'(t) = \frac{\int_{M_t}(H-h)(|A|^2 - H^2 +hH)d\mu}{\int_{M_t}d\mu}.
\ene
\el
\bp
An easy calculation using equations \eqref{mu} and \eqref{eq-H}. Note that the expression does not contain terms involving $\nabla H$.
\ep
The following inequalities for gradients are useful and we record them here.
\bl \label{HuiLemma1}(cf. \cite{Hui86}) The following inequalities hold: 
\begin{itemize}
\item[(i)] \beq |\nabla A|^2 \geq \frac{3}{n+2}|\nabla H|^2; \eeq
\item[(ii)] \beq |\nabla {\text\AA}|^2 \geq \frac{n-1}{2n+1}|\nabla A|^2 \geq \frac{3(n-1)}{(n+2)(2n+1)}|\nabla H|^2.\eeq
\end{itemize}
\el


\subsection{Intuitive decay of $|{\text\AA}|^2$} 

One of the key estimates for us is an exponential decay for $|{\text\AA}|^2$ on some {\ti}. We now give a heuristic argument to show why this is the 
case when $|{\text\AA}|^2$ is small and h-{\mcv}ity is preserved. 

Since $|{\text\AA}|^2 = |A|^2 - \frac{1}{n}H^2$ and $|\nabla {\text\AA}|^2 = |\nabla A|^2 - \frac{1}{n}|\nabla H|^2$, we obtain the {\ee} for $|{\text\AA}|^2$ as follows.
\bl \label{eq-tlA2-Cor}
\begin{align}\label{eq-tlA2}
\frac{\partial}{\partial t} |{\text\AA}|^2 
&= \Delta |{\text\AA}|^2 - 2|\nabla {\text\AA}|^2 + 2|{\text\AA}|^2(|A|^2 + n) + \frac{2h}{n}\bigl(H|A|^2 - n\,{\rm tr}(A^3)\bigr) \notag \\
&= \Delta |{\text\AA}|^2 - 2|\nabla {\text\AA}|^2 + 2|{\text\AA}|^2(|A|^2 + n) -2h\left\{ {\rm tr}({\text\AA}^3) + \frac 2n |{\text\AA}|^2H\right\} \,. 
\end{align}
\el
\bp
The {\ee} for $H^2$ can be easily derived from \eqref{eq-H}:
\be\label{eq-H2}
\frac{\partial}{\partial t} H^2 = \Delta H^2 -2|\nabla H|^2 + 2H(H-h)(|A|^2-n).
\ene
Then \eqref{eq-tlA2} follows easily from the identity (see e.g. page 335 of \cite{Li09}):
\beq
\text{tr}(A^3) -\frac1n |A|^2H = \text{tr}({\text\AA}^3) + \frac2n |{\text\AA}|^2H.
\eeq
\ep
To see a heuristic argument on exponential decay of $|{\text\AA}|^2$, we examine the equation \eqref{eq-tlA2} more closely, provided $|{\text\AA}|^2$ is small and $|h-H|$ is also very small. Obviously, one can apply the {\maxp} to \eqref{eq-tlA2} to obtain the exponential decay 
of $|{\text\AA}|^2$, if for some small $\epsilon>0$ we have  
$$2|{\text\AA}|^2(|A|^2 + n) -2h \left\{\text{tr}\left({\text\AA}^3\right) + \frac2n |{\text\AA}|^2H\right\}\leq -\epsilon|{\text\AA}|^2.$$ 
Since $|\text{tr}\left({\text\AA}^3\right)| \leq |{\text\AA}|^3$, it suffices to show 
\be\label{heu}
|{\text\AA}|^2 + \frac{H^2}{n} + n + |h|\cdot|{\text\AA}| - \frac{2hH}{n} < -\frac{\epsilon}{2}\,,
\ene
which can be rewritten as  
\beq
\frac{H^2}{n} =H+\frac{H(H-n)}{n}> |{\text\AA}|^2 + n + |h|\cdot |{\text\AA}| + \frac{2H(H-h)}{n}+\frac{\epsilon}{2}\,.
\eeq
This inequality holds once we establish $H > n+\sigma$ for some $\sigma>0$ (i.e., $h$-mean convexity) provided that $|{\text\AA}|^2$ and $|h-H|$ are both sufficiently small. We will make the argument precise in \S 3.3.


\subsection{Technical Tools}

For the sake of self-containedness of the paper, we now collect tools that will be used in the proof: a version of {\maxp}, Hamilton's interpolation inequalities for tensors, a generalization of 
Topping's theorem in {\hs}, and a $L^2$-bound for covariant derivatives of $A$ along the VPMCF. Firstly, the following version of {\maxp} is 
useful in our iteration scheme.
\bt\label{maxPrin} (Maximum Principle, see e.g. \cite[Lemma 2.11]{CLN06}) Suppose $u: M \times [0,T] \to \R$ satisfies
\beq 
\frac{\partial}{\partial t} u  \leq a^{ij}(t)\nabla_i\nabla_j u + \langle B(t), \nabla u\rangle + F(u)\,,
\eeq
where the coefficient matrix $\left(a^{ij}(t)\right)>0$ for all $t\in [0,T]$, $B(t)$ is a time-dependent vector field and $F$ is a Lipschitz function. If 
$u\leq c$ at $t=0$ for some $c\in\mathbb{R}$, then $u(x,t)\leq U(t)$ for all $(x,t)\in M \times \{t\}, t \in [0,T]$, where $U(t)$ is the solution to: 
\beq 
\frac{d}{dt} U(t) = F(U) \quad \text{with} \quad U(0) = c.
\eeq
\et

We also need the following Hamilton's interpolation inequalities for tensors. These inequalities will be used inductively for us to obtain integral 
bounds of covariant derives of ${\text\AA}$.
\bt(\cite{Ham82})\label{Hamiltonlemma2}
Let $M$ be an $n$-dimensional compact Riemannian manifold and $\Omega$ be any tensor on $M$. Suppose
\beq
\frac{1}{p} + \frac{1}{q} = \frac{1}{r} \quad \text{with } r\ge 1\,.
\eeq
Then we have the estimate:
\beq
\left(\int_{M} |\nabla \Omega|^{2r}\, d\mu\right)^{1/r} \leq (2r-2+n)\,
\left(\int_{M} |\nabla^2 \Omega|^p\, d\mu\right)^{1/p}\left(\int_{M} |\Omega|^{q}\, d\mu\right)^{1/q}\,.
\eeq
\et
\bt(\cite{Ham82})\label{Hamiltonlemma}
Let $M$ and $\Omega$ be the same as in the Theorem \ref{Hamiltonlemma2}. If $1\leq i\leq m-1$ and $m \ge 1$,
then there exists a constant $C = C(n,m)$ independent of the metric and connection on $M$, such that:
\beq
\int_{M} |\nabla^i \Omega|^{2m/i}\, d\mu \leq C\, \max_{M} |\Omega|^{2(m/i-1)}\int_{M}|\nabla^m \Omega|^2\, d\mu\,.
\eeq
\et

As an application of these inequalities and Corollary \ref{hdA-1}, we have the following:
\bl\label{interpolation01}
For any $m\geq 0$, we have the estimate
\begin{align*}
\left(\frac{d}{dt}\int_{M_t} |\nabla^m A|^2\right) &+ 2\int_{M_t} |\nabla^{m+1} A|^2
\leq C\,\max_{M_t}\left(1+|A| + |A|^2\right)\int_{M_t} |\nabla^m A|^2 \,,
\end{align*}
where $C=C(n,m, |h|)$.
\el
\bp
When $m=0$, the inequality is obvious in light of (\ref{eq-A2}). Now we consider $m\geq 1$. By integrating \eqref{hdA-2} of Corollary \ref{hdA-1} and 
using the generalized H\"{o}lder inequality we have:
\begin{align*}
&\left(\frac{d}{dt}\int_{M_t} |\nabla^m A|^2 \,d\mu\right) -\int_{M_t}(h-H)H|\nabla^m A|^2 \,d\mu + 2 \int_{M_t} |\nabla^{m+1} A|^2\,d\mu \\
\leq &C\,\Bigg\{\sum_{i+j+k=m}\left(\int_{M_t} |\nabla^i A|^{\frac{2m}{i}} \,d\mu\right)^{\frac{i}{2m}}\left(\int_{M_t} |\nabla^j A|^{\frac{2m}{j}} \,d\mu\right)^{\frac{j}{2m}} \left(\int_{M_t} |\nabla^k A|^{\frac{2m}{k}} \,d\mu\right)^{\frac{k}{2m}} \\
&\ \ + \sum_{r+s=m}\left(\int_{M_t} |\nabla^r A|^{\frac{2m}{r}} \,d\mu\right)^{\frac{r}{2m}}\left(\int_{M_t} |\nabla^s A|^{\frac{2m}{s}} \,d\mu\right)^{\frac{s}{2m}}\Bigg\}\left(\int_{M_t} |\nabla^m A|^2 \,d\mu\right)^{\frac12} \\
&\ \ \ \ + C\int_{M_t} |\nabla^m A|^2\,d\mu \,,
\end{align*}
where all the indices now take values from $1$ and up and the terms in the original sums with $0$ indices being absorbed by other sums and $C$'s.  

Applying Theorem \ref{Hamiltonlemma} for $A$, we have
\beq
\left(\int_{M_t} |\nabla^q A|^{2m/q} \,d\mu\right)^{q/2m} \leq C\, \max_{M_t}|A|^{1-q/m}\left(\int_{M_t} |\nabla^m A|^2 \,d\mu\right)^{1/2m}\,,
\eeq
where $q$ can be $i, j, k, r$ or $s$. We also notice  
\begin{align*}
\int_{M_t}|(h-H)H|\nabla^m A|^2 d\mu &\leq \max_{M_t}\{|h||H| + H^2\}\int_{M_t} |\nabla^m A|^2 d\mu   \\
& \leq C(n, |h|)\max_{M_t}\left(|A|^2 + |A|\right)\int_{M_t} |\nabla^m A|^2d\mu.
\end{align*}
Combining these inequalities, we complete the proof.
\ep

It's known from \cite{Hui87} that $L^\infty$ bound for $|A|$ along the VPMCF in Euclidean space and the initial $L^\infty$ bounds of its covariant derivatives will give 
$L^\infty$ bounds for the covariant derivatives. Adapting the argument there, we have the following lemma for explicit $L^2$ bounds for our situation.

\bl\label{le:L2}
Along the VPMCF, for $k\ge 0$, if 
$$\max\left\{|M_0|, \max_{M_t, t\in [0, T]}|A|^2\,, \,\max_{m\le k}\int_{M_0}|\nabla^m A|^2\,d\mu \right\}\,\leq \, \Lambda_0^2,$$ 
then uniformly for $t\in [0, T]$ and $m\le k$ we have
$$\int_{M_t}|\nabla^m A|^2\,d\mu\leq C(\Lambda_0, k)\,,$$
where $C(\Lambda_0, k)$ is independent of $T$. 
\el
     
\bp
Along the VPMCF, we have $|M_t|\le |M_0|$ by Corollary \ref{HuiLemma2}. So the conclusion is clear for $m=0$ for any fixed $k\geq 0$. We can then prove the lemma by induction on $m$. Suppose the conclusion is true for $m\geq 0$, to see this holds for $m+1 \leq k$, note that by Lemma \ref{interpolation01}, we know for $m\ge 0$,
\begin{align*}
\frac{d}{dt}\int_{M_t} |\nabla^m A|^2 \,d\mu
&\leq C(\Lambda_0)\int_{M_t} |\nabla^m A|^2 \,d\mu-2\int_{M_t} |\nabla^{m+1} A|^2\,d\mu\,,
\end{align*}
\begin{align*}
\frac{d}{dt}\int_{M_t} |\nabla^{m+1} A|^2 \,d\mu 
&\leq C(\Lambda_0)\int_{M_t} |\nabla^{m+1} A|^2 \,d\mu\,.
\end{align*}
Let $G(t)= C(\Lambda_0)\int_{M_t} |\nabla^m A|^2\,d\mu+\int_{M_t} |\nabla^{m+1} A|^2 \,d\mu$. Then we have 
\be \label{dddd}
G'(t) \leq C(\Lambda_0)\left(C(\Lambda_0)\int_{M_t} |\nabla^m A|^2 \,d\mu-\int_{M_t} |\nabla^{m+1} A|^2\,d\mu\right)\,.
\ene
Consider the maximum of $G(t)$ achieved at $t=\bar{t} \in [0,T]$. If $\bar{t} =0$ then for all $t\in [0,T]$,
\be \label{dddd1}
G(t) \leq G(0) \leq (C(\Lambda_0)+1)\Lambda_0^2\,.
\ene
Otherwise by \eqref{dddd},
$$C(\Lambda_0)\int_{M_{\bar{t}}} |\nabla^m A|^2 \,d\mu-\int_{M_{\bar{t}}} |\nabla^{m+1} A|^2\,d\mu\geq 0\,,
$$
and thus we have  for all $t\in [0,T]$ 
\be \label{dddd2}
G(t) \leq G(\bar{t}) \leq C(\Lambda_0, k)\,.
\ene
Therefore by \eqref{dddd1} and \eqref{dddd2},
$$
\int_{M_t} |\nabla^{m+1} A|^2\,d\mu \leq C(\Lambda_0, k),
$$
which is independent of $T$.
\ep

In the Euclidean space, Topping (\cite{Top08}) discovered a relation between the intrinsic diameter and the {\mc} $H$ of any closed, 
connected and smoothly immersed submanifold. This result has been extended to a more general Riemannian setting by 
Wu-Zheng (\cite{WZ11}), using Hoffman-Spruck's generalization (\cite{HS74}) of the Michael-Simon's inequality (\cite{MS73}). 
We formulate their result in our setting below. 
\bt(\cite{WZ11})\label{top}
Let $M$ be an $n$-dimensional closed, connected manifold smoothly isometrically immersed in
$\H^{N}$, where $N \geq n+1$. There exists a constant $C=C(n)$ such that the 
intrinsic diameter and the {\mc} $H$ of $M$ are related by the following inequality:
\beq
{\rm diam}\,(M) \leq C(n)\int_M |H|^{n-1}\,d\mu\,.
\eeq
\et
\subsection{Hyperbolic mean convexity}

The $h$-{\mcv}ity is a very important geometric ingredient in our main result. Note that mean convexity and $h$-{\mcv}ity are not known to be preserved 
along the VPMCF. The strict $h$-convexity is however preserved along the VPMCF in $\H^{n+1}$ \cite{CRM07}. We give an alternative proof for this result 
by following very closely Huisken's tensor calculations in \cite{Hui84, Hui87} and highlighting the role of the curvature for the ambient space. Unlike the preserved mean convexity along the MCF in Euclidean space, this shows the subtlety of the $h$-{\mcv}ity in hyperbolic space and the negative-curvature effects of the ambient space.

\begin{pro}\label{hconvex} 
(\cite{CRM07}) Let $M^n$ be a smooth, embedded, closed hypersurface moving by the VPMCF \eqref{eq-vp} in a smooth, complete, hyperbolic manifold 
$N^{n+1}$. If the initial hypersurface $M^n$ is strictly $h$-convex, then each evolving hypersurface $M_t^n$ is also strictly h-convex along the flow \eqref{eq-vp}.
\end{pro}

\bp Let $M_{ij} = a_{ij} - g_{ij}$. Recall the {\ee}s for $a_{ij}$ and $g_{ij}$ along the {\mcf} \eqref{eq-vp} as (\ref{eq-a}) and (\ref{eq-g}):
\beqar
\ppl{}{t}a_{ij} - \Delta a_{ij} &=& (h-2H)a_{i\ell}a_{j\ell} +|A|^2a_{ij} -n a_{ij} - h\bar{R}_{0i0j} \\
&&\ \ \ \ \ \ \ -a_{j\ell}\bar{R}_{\ell mim}-a_{i\ell}\bar{R}_{\ell mjm}+2a_{\ell m}\bar{R}_{\ell imj},
\eeqar
where the covariant derivatives for the curvature tensor disappear since the sectional curvature is $-1$, and 
\beq
\ppl{}{t}g_{ij} = 2(h-H)a_{ij}.
\eeq
Therefore we obtain the {\ee} for the symmetric tensor $M_{ij}$:
\beq
\ppl{}{t}M_{ij} = \Delta M_{ij} + N_{ij},
\eeq
where we have used $\Delta g=0$ and  
\bear\label{Nij}
N_{ij} &=& (h-2H)a_{i\ell}a_{j\ell} +(|A|^2-n)a_{ij}- h\bar{R}_{0i0j}-a_{j\ell}\bar{R}_{\ell mim}\notag\\
&&\ \ \ \ \ \ \ -a_{i\ell}\bar{R}_{\ell mjm}+2a_{\ell m}\bar{R}_{\ell imj} + 2(H-h)a_{ij}.
\eear
Now recall from \eqref{rie-tensor},  
\beq
\bar{R}_{\alpha\beta\gamma\delta} = (-1)\cdot(\bar{g}_{\alpha\gamma}\bar{g}_{\beta\delta}-
\bar{g}_{\alpha\delta}\bar{g}_{\gamma\beta}).
\eeq
Let $X$ be a null-eigenvector of $M_{ij}$ at some $(x_0,t_0)$. We arrange the coordinates such that at 
$(x_0,t_0)$, $X= e_1$, $g_{ij} = \delta_{ij}$ and $a_{ij} = \lambda_i\delta_{ij}$. This is justified as 
$\{g_{ij}\}$ is a symmetric positive-definite matrix, $\{a_{ij}\}$ is a symmetric matrix, and so they  
can be simultaneously diagonalized.

We examine term by term from \eqref{Nij} to arrive at:
\beq
N_{11} =(h-2H)\lambda_1^2+(|A|^2-n)\lambda_1 +h+ 2(n-1)\lambda_1 +2(\lambda_1-H)+2(H-h)\lambda_1.
\eeq
Meanwhile, with $X=e_1$ being a null-eigenvector of $M_{ij}$, we have $\lambda_1 = 1$ since $M_{11}=a_{11} - g_{11}=0$. Thus, we have  
\beq
 N_{11} =|A|^2+n-2H \ge \frac1n H^2 -2H + n = \frac1n (H-n)^2 \ge 0.
\eeq
The conclusion follows from Hamilton's {\maxp} for tensors (\cite{Ham82}).
\ep


\section{Proof of Main Theorem}

We are now ready to use iteration method to prove our main theorem. It's divided into four steps discussed in four  subsections accordingly.


\subsection{Step One: Short Time Bounds} 

We start by bounding important geometric quantities for short time, with the bounds depending on the initial conditions. This is certainly 
expected for a smooth flow. However, one expects such bounds to hold only for a short time, and as the flow evolves such bounds would deteriorate by extending the time interval.

The first technical lemma is as follow:
 
\bl\label{key}
Let $M_t^n \subset \H^{n+1} $, $n\geq 2$ be a smooth closed solution to the VPMCF \eqref{eq-vp} for $t\in [0,T)$ with $T\leq \infty$. Assume 
\be\label{Condition-01}
\max\left\{|M_0|^2\,, \, \max_{M_0}|A|^2\,, \,\int_{M_0}|\nabla^m A|^2\,d\mu \right\}\,\leq \, \Lambda_0^2
\ene
for some $\Lambda_0 \gg 1$ and all $m \in [1, n+3]$, where $|M_t|$ is the 
$n$-dimensional surface area of $M_t$ with the induced metric. There exist constants
$\epsilon_0 = \epsilon_0(n,\Lambda_0)>0$ and $t_1 = t_1(n,\Lambda_0) \in (0,1)$ such that if
\be\label{Condition-02}
\int_{M_0}|{\text\AA}|^2 \,d\mu \leq \epsilon < \epsilon_0\,,
\ene
then for any $t\in [0,t_1]$ and any $m\in [0,n+3]$ we have
\be\label{estimate-1}
\max\left\{\max_{M_t}|A|^2,\,\int_{M_t}|\nabla^m A|^2\,d\mu \right\}\leq  2\Lambda_0^2\,.
\ene
Moreover, there exist $C_1 = C_1 (n, \Lambda_0)$ and some universal constant $\alpha \in (0,1)$ such that 
for any $t\in [0,t_1]$
\be\label{estimate-2}
\max_{M_t}\left(|{\text\AA}|+ |\nabla H| + |h- H|\right) \leq C_1 \epsilon^{\alpha}\,.
\ene
\el

\bp
Recall from \eqref{eq-A2} the {\ee} for $|A|^2$ is given by
\beq
\frac{\partial}{\partial t} |A|^2 = \Delta |A|^2 - 2|\nabla A|^2 + 2|A|^2(|A|^2+n) - 2h\,\text{tr}\left(A^3\right) + 2H(h-2H).
\eeq
Using the facts that $|\text{tr}\left(A^3\right)|\leq |A|^3$ (see Lemma 2.2 \cite{HS99}), and $H^2 \leq n |A|^2$, we 
obtain the following inequality on $M_t$ for all $t \in [0,T)$:
\be
\frac{\partial}{\partial t} |A|^2 \leq \Delta |A|^2 + 2|A|^4 + 2n|A|^2 + 2|h|(|A|^3+\sqrt{n}|A|).
\ene
Set $f(t) = \max\limits_{M_t}|A|^2$, then $f(t)$ satisfies 
\bear
\frac{\partial}{\partial t} f &\leq& 2f^2+ 2nf + 2|h|(|A|^3+\sqrt{n}|A|) \nonumber \\
&\leq& 2f^2+ 2nf + 2\sqrt{n}f^2+ 2nf \nonumber \\
&\leq& 4nf^2+ 4nf. 
\eear

One solves the comparison ODE explicitly to get $U(t) >0$ satisfying
\beq
\log\left(1+\frac{1}{U(t)}\right) = \log\left(1+\frac{1}{U(0)}\right) - 4nt,
\eeq
with $U(0)= f(0)=\max\limits_{M_0} |A|^2\leq \Lambda_0^2$ by (\ref{Condition-01}). So $f(t) \leq U(t)$ for all $t \in [0,T)$. 

Therefore, there exists some $ t_1 = t_1 (n,\Lambda_0) \in (0,1)$ such that
\be\label{Hsqbound-2}
\max_{M_t}|A|^2 \leq 2\Lambda_0^2\quad \text{for all } t\in [0, t_1]\,.
\ene
Moreover, by choosing $t_1$ sufficiently small and integrating the inequality in Lemma \ref{interpolation01} 
over $[0,t_1]$, we have 
\be\label{covariantbounds}
\int_{M_t} |\nabla^m A|^2 \,d\mu \leq e^{C(n, \Lambda_0)t_1} \int_{M_0} |\nabla^m A|^2 \,d\mu \leq 2\Lambda_0^2
\ene
for all $t\in [0,t_1]$ and $m \in [1, n+3]$. Using the Sobolev embedding on compact manifolds \cite{Aub98}, this yields 
\be\label{gradientbound1}
|A|_{C^2(M_t)}  \leq C(n, \Lambda_0) \quad \text{for all }\, t\in [0,t_1]\,.
\ene

In light of 
\beq
|h|\leq \max_{M_t}|H|\leq \sqrt{n}\max_{M_t}|A| \leq \sqrt{2n}\Lambda_0,
\eeq 
\beq 
|\text{tr}({\text\AA}^3)|\leq |{\text\AA}|^3 \leq \sqrt{2}\Lambda_0|{\text\AA}|^2,
\eeq 
we integrate the {\ee} \eqref{eq-tlA2} for $|{\text\AA}|^2$ over $M_t$ for $t\in [0,t_1]$ to get
\be\label{good-1}
\frac{\partial}{\partial t} \int_{M_t}|{\text\AA}|^2 \, d\mu \leq C(n, \Lambda_0)\int_{M_t}|{\text\AA}|^2 \, d\mu\,,
\ene
and so using (\ref{Condition-02}) we have
\be\label{good2}
 \int_{M_t}|{\text\AA}|^2 \, d\mu \leq \epsilon e^{C(n, \Lambda_0)t} \leq C(n, \Lambda_0)\epsilon \quad \text{for all}~t\in [0,t_1]\,,
\ene
where the constant $C(n, \Lambda_0)$ can be different at places. We then apply Hamilton's interpolation inequalities 
(Theorem \ref{Hamiltonlemma2} with $r=1,\ p=q=2$):
\be\label{good-eq-6}
\int_{M_t} |\nabla {\text\AA}|^2\,d\mu \leq n \left(\int_{M_t} |{\text\AA}|^2\,d\mu\right)^{\frac{1}{2}}
\left(\int_{M_t} |\nabla^2 {\text\AA}|^2\,d\mu\right)^{\frac{1}{2}} \leq C(n,\Lambda_0)\epsilon^{\frac{1}{2}}\,,
\ene
where we use $|\nabla^2 {\text\AA}|\leq C(n)|\nabla^2 A|$ and the $L^2$-bound for $|\nabla^2 A|$ in \eqref{covariantbounds}. In fact, applying Theorem \ref{Hamiltonlemma2} inductively, 
we have for all $m \in [0,n+2]$,
\be
\int_{M_t} |\nabla^m {\text\AA}|^2\,d\mu \leq C(n, \Lambda_0)\epsilon^{1/2^m} \quad \text{for all } t\in [0,t_1]\,.
\ene
Now again by the Sobolev embedding \cite{Aub98}, we have:
\be\label{good3}
|{\text\AA}|_{C^2(M_t)} \leq C(n,\Lambda_0)\epsilon^{\alpha},
\ene
for all $t\in [0,t_1]$ and some universal constant $\alpha \in (0,1)$. Now by (ii) of Lemma \ref{HuiLemma1}, for all $t\in [0,t_1]$ we have
\be\label{Hgrad1}
\max_{M_t}|\nabla H| \leq C(n)\max_{M_t}|\nabla {\text\AA}| \leq C(n, \Lambda_0)\epsilon^{\alpha} \,.
\ene
Furthermore, by Corollary \ref{HuiLemma2}, the surface area $|M_t|$ is non-increasing along the flow, i.e.
\be\label{Hsqbound-30}
|M_t| \leq |M_0| \leq \Lambda_0^2\,.
\ene
Using Theorem \ref{top}, \eqref{Hsqbound-2}, \eqref{Hgrad1} and \eqref{Hsqbound-30}, we arrive at
\begin{align}\label{hHbound02}
|h(t)-H(x,t)| &= \left(\int_{M_t}\, d\mu\right)^{-1} \left|\int_{M_t} H(y,t) \,  - H(x,t) d\mu(y)\right|\notag\\
& \leq \text{diam}\,(M_t)\max_{M_t} |\nabla H|\\
&\leq C(n, \Lambda_0)\epsilon^{\alpha}\notag
\end{align}
for all $(x,t)\in M_t$ and $t\in [0,t_1]$. This together with \eqref{good3} and \eqref{Hgrad1} give \eqref{estimate-2}, and we conclude the proof.
\ep

With the above control of geometric quantities, we next show that the $h$-{\mcv}ity is preserved for short time if the initial hypersurface is close to an umbilical sphere in the $L^2$-sense. 

\bl \label{key2}
Let $M_t^n \subset \H^{n+1}$ for $n\geq 2$ be a smooth closed solution to the VPMCF \eqref{eq-vp} as in 
Lemma \ref{key} with the initial condition \eqref{Condition-01}. Suppose 
\be \label{mean-convex-1}
\min_{M_0} (H - n) \ge c_0 > 0\,.
\ene
Then there exist $\epsilon_1 = \epsilon_1(n,\Lambda_0) \in (0, \epsilon_0)$ and $T_1 = T_1(n,\Lambda_0) \in (0,t_1]$, where $\epsilon_0$ and $t_1$ are as in Lemma \ref{key}, such that if
\beq
\int_{M_0}|{\text\AA}|^2 \,d\mu \leq \epsilon < \epsilon_1\,,
\eeq
then for $t\in [0,T_1]$ we have
\be\label{estimate-3}
\min_{M_t} (H - n) \ge \frac{c_0}{2} > 0\,.
\ene
\el
\bp
We start with the {\ee} for $H$ \eqref{eq-H}:
\beq
H_t = \Delta H + (H-h)(|A|^2-n).
\eeq
By \eqref{Hsqbound-2} and \eqref{gradientbound1}, for any $(x,t)\in M_t$, $t\in [0,t_1]$, we have:
\be\label{Hgrad2}
\left|\frac{\partial}{\partial t} H\right|(x,t) \leq C(n, \Lambda_0),
\ene
where we have also used $|\nabla^2 H|\leq C(n)|\nabla^2 A|$. Using \eqref{Hgrad1} and \eqref{Hgrad2} and choosing $T_1 = T_1(n, \Lambda_0) \in (0,t_1]$ and $\epsilon_1 = \epsilon_1(n,\Lambda_0) \in (0, \epsilon_0)$ sufficiently small, we have
\beq
\min_{M_t} (H - n) \geq \frac{1}{2}\min_{M_0} (H - n) \geq \frac{c_0}{2}> 0\,.
\eeq
\ep


\subsection{Step Two: Reduction}

In the previous subsection we have obtained estimates \eqref{estimate-1} and \eqref{estimate-2} on some {\ti} $[0,t_1]$, provided that the initial hypersurface is close to an umbilical sphere in the $L^2$-sense (see \eqref{Condition-2}). In this step, we make a key reduction. Namely, we show it suffices to prove the main theorem when the mean curvature $H$ of the evolving hypersurface is close to $n$. In particular, we have the following. 

\begin{pro}\label{red}
Let $M_t^n \subset \H^{n+1}$ for $n \geq 2$ be a smooth closed solution to the VPMCF \eqref{eq-vp} on 
$t\in [0,t_1]$ with $t_1= t_1(n,\Lambda_0)\in (0,1)$, where $t_1$ and $\Lambda_0$ are as in Lemma ~\ref{key}. If \eqref{Condition-01} and \eqref{Condition-02} hold, then 
\begin{enumerate}

\item either the evolving hypersurface $M_t$ becomes strictly h-convex, and the flow \eqref{eq-vp} exists for 
all time and converges exponentially to an umbilical sphere, 
 
\item or there is a constant $C_2 = C_2 (n, \Lambda_0) > 0$ such that for all $(x,t)\in M_t$, $t\in [0,t_1]$ we have 
\be\label{H-n}
|H(x, t) -n| \le C_2\epsilon^{\frac{\alpha}{2}},
\ene
where $\epsilon$ is from \eqref{Condition-02} and $\alpha \in (0,1)$ is from \eqref{estimate-2}.
\end{enumerate}
\end{pro}

\bp
On the {\ti} $[0,t_1]$, we recall the estimate \eqref{estimate-2} from Lemma ~\ref{key}:
\beq
\max_{M_t}\left(|{\text\AA}|+ |\nabla H| + |h- H|\right) \leq C_1 \epsilon^{\alpha}
\eeq 
for some $C_1 = C_1(n, \Lambda_0) > 0$. Let $\{\lambda_i\}_{i = 1, 2, \cdots, n}$ be the {\pc}s of $M_t$ at $(x,t)\in M_t$. Direct algebra gives  
\be\label{pc}
|{\text\AA}|^2 = \frac1n \sum_{i<j}(\lambda_i - \lambda_j)^2, 
\ene
so there exists $C_3 = C_3(n,\Lambda_0)>0$ such that for all $(x,t) \in M_t$, $t\in [0,t_1]$,
\be\label{pc-2}
|\lambda_i (x,t)- \lambda_j(x,t)| \le C_3\epsilon^{\alpha}.
\ene
Therefore for all $(x,t) \in M_t$, $t\in [0,t_1]$ and any fixed $i \in \{1,2, \cdots, n\}$, we have 
\be\label{H-n2}
|H(x,t) - n\lambda_i(x)| \le C_4\epsilon^{\alpha},
\ene
for some $C_4 = C_4(n,\Lambda_0) > 0$.

For some $C_5 = C_5(n,\Lambda_0) > 0$ which will be fixed shortly, suppose there is $\eta_0=C_5\epsilon^{\frac{\alpha}{2}}> 0$ where $\epsilon \in (0, \epsilon_0)$ 
and some $(x_0,t_0) \in M_{t_0}$ where $t_0\in [0,t_1]$ such that $H(x_0,t_0) < n -\eta_0$. Then from \eqref{H-n2} we have:
\beq
n\lambda_i(x_0,t_0) - C_4\epsilon^\alpha \le H(x_0,t_0) < n-\eta_0=n-C_5\epsilon^{\frac{\alpha}{2}}.
\eeq 
Since $\epsilon \in (0,\epsilon_0)$ is small, for properly chosen $C_5$  and $C_6=C_6(n, \Lambda_0)>0$, we have $\lambda_i (x_0, t_0)< 1-C_6\epsilon^{\frac\alpha 2}$ for all $i \in \{1,2, \cdots, n\}$. In light of $\max_{M_t}|\nabla H| \leq C_1 \epsilon^{\alpha}$, the smallness of $\epsilon$ and the diameter bound from Theorem \ref{top}, we have $H < n$ at every point of $M_{t_0}$. However, this contradicts the fact that any smooth closed 
hypersurface has at least one point whose {\mc} is greater than $n$ in $\H^{n+1}$ by comparing with horospheres. 

Similarly for some $C'_5 = C'_5(n,\Lambda_0)>0$ which will be fixed shortly, suppose there is some $\eta'_0 =C'_5\epsilon^{\frac{\alpha}{2}}> 0$ where $\epsilon \in (0, \epsilon_0)$ and some $(x'_0, t'_0) \in M_{t'_0}$ such that $H(x'_0, t'_0) > n +\eta'_0$. We have
\beq
n\lambda_i(x'_0, t'_0) + C_4\epsilon^\alpha \ge H(x'_0, t'_0) > n +\eta'_0=n+C'_5\epsilon^{\frac{\alpha}{2}}.
\eeq 
Using again the smallness of $\epsilon$, for properly chosen $C'_5$ and $C'_6$ we have $\lambda_i (x'_0, t'_0)> 1+C'_6\epsilon^{\frac\alpha 2}$ for any $i \in \{1,2, \cdots, n\}$. Using again the fact that $\max_{M_t}|\nabla H| \leq C_1 \epsilon^{\alpha}$, smallness of $\epsilon$ and the diameter bound from Theorem \ref{top}, we find $\lambda_i (x,t'_0)> 1$ for all $i \in \{1,2, \cdots, n\}$ and all $(x,t'_0) \in M_{t'_0}$. Namely, $M_{t'_0}$ is strictly {\it h-convex}. By the main theorem of \cite{CRM07}, the VPMCF then exists for all time after $t=t'_0$, stays strictly {\it h-convex} and converges exponentially to an umbilical sphere in $\H^{n+1}$.

Finally, we are left with \eqref{H-n}, which completes the proof. 
\ep

\br
By Proposition ~\ref{red}, we can now assume $H$ of $M_t$ is very close to $n$ on {\ti} $[0,t_1]$, namely the inequality \eqref{H-n}, for the remaining proof for Theorem \ref{main}, and therefore we now have $H >0$ (hence $h>0$). 
\er


\subsection{Step Three: Precise Decay} In the previous subsection we have obtained estimates \eqref{estimate-1}, \eqref{estimate-2} and \eqref{estimate-3} 
on some short time interval $[0,T_1]$, provided that the initial hypersurface is close to an umbilical sphere in the $L^2$ sense (see \eqref{Condition-2}) and $h$-mean convex (see \eqref{mean-convex-1}). These bounds 
will likely deteriorate along the flow if we iterate for later time intervals. For an iteration argument to work, we need to establish time-independent 
bound on these quantities for this short time interval. 

In this subsection, we show that, if estimates similar to \eqref{estimate-1}, 
\eqref{estimate-2} and \eqref{estimate-3} hold on some time interval $[0,T_1]$, then we can choose sufficiently small $\epsilon$ in the initial $L^2$-bound \eqref{Condition-2} on ${\text\AA}$, such that $|{\text\AA}|$, $|\nabla H|$ and $|h-H|$ exponentially decay on this time interval $[0,T_1]$. More precisely, we establish the following theorem.

\bt\label{exp1}
Let $M_t^n \subset \mathbb{H}^{n+1}$ for $n\geq 2$ be a smooth closed solution to the VPMCF \eqref{eq-vp} 
with the initial condition
\beq
\int_{M_0}|{\text\AA}|^2 \,d\mu \leq \epsilon\,.
\eeq
Suppose for any $t\in [0,T_1]$ with $T_1\leq \infty$ and all $m\in [1, n+3]$ we have
\be\label{Condition-3}
\max\left\{|M_0|^2\,, \, \max_{M_t}|A|^2\,, \,\int_{M_0}|\nabla^m A|^2\,d\mu \right\} \leq \Lambda_1^2, ~~~~\min_{M_t} (H - n)\geq \sigma,
\ene
\be\label{Condition-4}
\max_{M_t}\left(|{\text\AA}| + |\nabla H| + |h- H|\right) \leq C_1 \epsilon^{\beta},
\ene
for constants $\Lambda_1>0, \sigma>0, \beta\in (0, 1)$ and $C_1>0$. Then there exists some $\epsilon_2 = \epsilon_2 (n, \Lambda_0,\beta, C_1)>0$ such that if $\epsilon <\epsilon_2$, then for all $t\in [0,T_1]$ we have 
\be\label{TracelessADecay1}
\max_{M_t}|{\text\AA}| \leq \max_{M_0}|{\text\AA}|, 
\ene
\be\label{hHDecay1}
\max_{M_t} \left(|{\text\AA}| + |\nabla H| + | h- H|\right) \leq
C_2(n,\Lambda_1, C_1)\left(\max_{M_0}|{\text\AA}|\right)^{\alpha} e^{-\alpha\sigma t}\,,
\ene
where $\alpha\in (0, 1)$ is the universal constant from Lemma \ref{key}. 
\et

\bp

To start with, by Lemma \ref{le:L2} and (\ref{Condition-3}), for $m\in [1, n+3]$ and $t\in [0,T_1]$ we have  
$$\int_{M_t}|\nabla^m A|^2\,d\mu\leq C(n, \Lambda_1),$$
which works as the replacement of (\ref{estimate-1}) as in the proof of Lemma \ref{key}. Now using \eqref{Condition-3} we compute
\bear\label{estimate-4}
n - \frac{hH}{n} &=& n - \frac{H\int_{M_t} H \, d\mu}{n \int_{M_t} \, d\mu} \nonumber\\
&\le& n-\frac{(n+\sigma)^2}{n}  \nonumber\\
&<& -2\sigma\,,
\eear
and 
\begin{align}\label{estimate-5}
\left|\frac{1}{n}H^2 - \frac{hH}{n}\right|(x,t) &= \left|H(x,t) \cdot \frac{\int_{M_t}\left[ H(x,t)- H(y,t)\right]\, d\mu(y)}{n\int_{M_t}\, d\mu}\right| \notag\\
&\leq \frac{1}{n}\max_{M_t}H \cdot \text{diam}\,(M_t)\cdot\max_{M_t}|\nabla H| \notag\\
&\leq C(n, \Lambda_1, C_1)\epsilon^\beta\,,
\end{align}
where we have used $|H| \leq \sqrt{n} |A| \leq \sqrt{n} \Lambda_1$ and Theorem \ref{top}.

Now by \eqref{eq-tlA2}, \eqref{estimate-4} and \eqref{estimate-5}, we have
\begin{align*}
\frac{\partial}{\partial t} |{\text\AA}|^2 &= \Delta |{\text\AA}|^2 - 2|\nabla {\text\AA}|^2+2|{\text\AA}|^2(|A|^2 + n) -2h \left\{\text{tr}\left({\text\AA}^3\right) + \frac2n |{\text\AA}|^2H\right\}\\
&\leq \Delta |{\text\AA}|^2+ 2|{\text\AA}|^2(|{\text\AA}|^2 + \frac{1}{n}H^2 + n) + 2h \left|{\text\AA}\right|^3 - \frac{4hH}{n} |{\text\AA}|^2\\
&= \Delta |{\text\AA}|^2 + 2\left(|{\text\AA}|^2 + h |{\text\AA}| + \frac{1}{n}H^2 + n - \frac{2hH}{n}\right)|{\text\AA}|^2\\
&\leq \Delta |{\text\AA}|^2 - (4\sigma-\widehat C\epsilon^\beta)|{\text\AA}|^2\\
&\leq \Delta |{\text\AA}|^2 - \sigma|{\text\AA}|^2. 
\end{align*}
where $\widehat C = \widehat C(n, \Lambda_1, C_1)>0$ and for the the last step we choose $\epsilon$ to be sufficiently small. Therefore, we conclude the 
exponential decay of $|{\text\AA}|$ from the {\maxp}, i.e. Theorem \ref{maxPrin}, $$\max_{M_t}|{\text\AA}|^2\leq e^{-\sigma t}\max_{M_0}|{\text\AA}|^2,$$
and the estimate \eqref{TracelessADecay1} also follows. This is where the $h$-mean convexity is essentially involved in our arguments, see \eqref{estimate-4}. Afterwards, we can prove \eqref{hHDecay1} by the exact arguments in the proof of 
Lemma ~\ref{key}, namely \eqref{good-eq-6}--\eqref{hHbound02}.
\ep


\subsection{Step Four: Time Extension}

In this step, we use the exponential decay of $|{\text\AA}|$, $|\nabla H|$ and $|h-H|$ on some short time interval 
obtained in previous step to extend the {\ti} of interest.

\bt\label{ext}
Let $M_t^n \subset \mathbb{H}^{n+1}$ for $n\geq 2$ be a smooth closed solution to the VPMCF \eqref{eq-vp} with the initial hypersurface satisfying 
$$|M_0|\leq \Lambda_0, ~~\max_{M_0} |H| \leq \Lambda_0, ~~\int_{M_0}|\nabla^m A|^2\,d\mu\leq \Lambda_0^2, ~~\min_{M_0} (H - n)\geq \frac{1}{\Lambda_0^2} >0$$
for all $m\in [1,n+3]$. Suppose for any $t\in [0,T]$ with $T<\infty$ we have 
\be\label{wholebound-1}
\max_{M_t}|A|^2 \leq \, \Lambda_0^2, ~~~~\min_{M_t} (H - n) \geq \frac{1}{2\Lambda_0^2} >0
\ene
and
\be\label{wholebound-2}
\max_{M_t}\left(|{\text\AA}| + |\nabla H| + | h- H|\right) \leq C_\ast \epsilon^\frac{\alpha^2}{2} e^{-\alpha\sigma t} \leq C_\ast \epsilon^\frac{\alpha^2}{2}\, ,
\ene
where $\alpha\in(0, 1)$ is the universal constant from Lemma ~\ref{key} and $\sigma=\frac{1}{2\Lambda_0^2}$ is as in Theorem \ref{exp1}. Then there exist $\epsilon_3 = \epsilon_3(n, \Lambda_0, \alpha, C_\ast)>0$ and $T_2 = T_2(n,\Lambda_0) >0$ such that if
\be
\int_{M_0}|{\text\AA}|^2 \,d\mu \leq \epsilon < \epsilon_3\,,
\ene
then \eqref{wholebound-1} and \eqref{wholebound-2} hold for $t\in [0,T+T_2]$.
\et

\bp 
We begin by applying Lemma \ref{key} and Lemma \ref{key2} to obtain $\epsilon_4= \epsilon_0(n,\Lambda_0^2)$ and $T_2= T_1(n,\Lambda_0^2)$ such that if
\beq
\int_{M_0}|{\text\AA}|^2 \,d\mu \leq \epsilon <\epsilon_4\,,
\eeq
then for all $t\in [T,T+T_2]$ we have
\be\label{good5}
\max\left\{\max_{M_t}|A|^2\,, \,\int_{M_t}|\nabla^m A|^2\,d\mu \right\} \leq  2\Lambda_0^2 \quad \text{and}\quad
\min_{M_t} (H - n) \geq \frac{1}{4\Lambda_0^2},
\ene
\be\label{good6}
\max_{M_t}\left(|{\text\AA}|+ |\nabla H| + |h- H|\right) \leq C_1(n,\Lambda_0) \epsilon^{\alpha}\,,
\ene
where $C_1$ and $\alpha$ are from Lemma \ref{key}. Then choose $\epsilon_5 = \epsilon_5(n,\Lambda_0, \alpha, C_\ast)>0$ sufficiently small so that for any $\epsilon<\epsilon_5$, we have
\beq
C_1(n,\Lambda_0) \epsilon^{\alpha - \frac{\alpha^2}{2}}\leq C_\ast.
\eeq
Therefore for all $t\in[0,T+T_2]$ we have \eqref{good5} and also
\be\label{good7}
\max_{M_t}\left(|{\text\AA}| + |\nabla H| + | h- H|\right)  \leq C_\ast \epsilon^\frac{\alpha^2}{2}\,.
\ene

By Corollary \ref{HuiLemma2} the surface area $|M_t|$ is non-increasing along the flow, therefore 
$|M_t|\leq \Lambda_0 <\Lambda_0^2$ by the initial condition \eqref{Condition-1} as long as the flow exists, 
in particular, on $[0,T+T_2]$. Now we apply the Theorem ~\ref{exp1} on $[0, T+T_2]$ with 
$\Lambda_1^2 = 2\Lambda_0^2$, $C_1 = C_\ast$, $\beta = \frac{\alpha^2}{2}$ and $\sigma=\frac{1}{4\Lambda_0^2}$ to conclude that for some $\epsilon_6:= \epsilon_2 (n,\Lambda_0,\alpha, C_\ast)>0$ sufficiently small, if $\epsilon < \epsilon_6$, then for all $t\in [0,T+T_2]$, we have 
\begin{align}
\max_{M_t} \left(|{\text\AA}| + |\nabla H| + | h- H|\right) 
&\leq C_2(n, \Lambda_0, C_\ast)\left(\max_{M_0}|{\text\AA}|\right)^{\alpha}e^{-\alpha\sigma t}\notag \\
&\leq C_2(n, \Lambda_0, C_\ast)[C_1(n, \Lambda_0) \epsilon^{\alpha}]^{\alpha}e^{-\alpha\sigma t}\,,
\end{align}
where we've used \eqref{estimate-2}  at $t=0$. Now choose $\epsilon_7 = \epsilon_7 (n, \Lambda_0, \alpha, C_\ast)>0$ small enough so that
\be\label{good8}
 C_2(n,2\Lambda_0^2,C_\ast)[C_1(n, \Lambda_0)]^{\alpha} \epsilon^{\frac{{\alpha}^2}{2}} \leq C_\ast, 
\ene
thus \eqref{wholebound-2} holds for all $t\in [0,T + T_2]$. 

We are left to show \eqref{wholebound-1} for $t\in [0,T + T_2]$. Let's examine each term in \eqref{wholebound-1}. Consider $\max_{M_t}|A|$. Recall the time derivative formula for $h(t)$ \eqref{derivativeh} is given by
\beq
h'(t) = \frac{\int_{M_t}(H-h)(|A|^2 - H^2 +hH)d\mu}{\int_{M_t}d\mu}.
\eeq
Then using \eqref{good5} and \eqref{good6}, we have   
\beq
|h'(t)| \leq C_3(n, C_\ast, \Lambda_0) \epsilon^\frac{\alpha^2}{2}e^{-\alpha\sigma t}
\eeq
for all $t\in [0,T+T_2]$. Note that, from the initial condition \eqref{Condition-1}, we also have
\beq
h(0) = \frac{\int_{M_0} H\, d\mu}{\int_{M_0} \, d\mu} \leq \max_{M_0}|H| \leq \Lambda_0\,.
\eeq
By choosing $\epsilon < \epsilon_8 = \epsilon_8(n, \Lambda_0, \alpha, C_\ast)$ sufficiently small, we then have for any $t\in [0,T+T_2]$:
\be
|h(t)| \leq \frac65\Lambda_0\,.
\ene
Then by \eqref{good6} and $n\geq 2$, for sufficiently large $\Lambda_0$ we have
\be
\max_{M_t} |A| = \max_{M_t} \sqrt{|{\text\AA}|^2 + \frac{1}{n}H^2} \leq 
\max_{M_t}\left(|{\text\AA}|+ \frac{1}{\sqrt{n}}|H-h|\right)+ \frac{1}{\sqrt{n}}|h(t)|\leq \Lambda_0 \,.
\ene

Finally, we consider the term $\min_{M_t} (H - n)$. Using the {\ee}s for $H$ (see \eqref{eq-H}) and $d\mu$ (see \eqref{mu}), we have
\begin{align}
\int_{M_t} H\, d\mu - \int_{M_0} H\, d\mu &= \int_0^t\int_{M_s} H^2(h-H) + (H-h)(|A|^2-n)\,d\mu \,ds \notag \\
&\geq - C(n, \Lambda_0, C_\ast) \epsilon^\frac{\alpha^2}{2}\int_0^t  e^{-\alpha\sigma s}\,ds \geq - 
C_4(n, \Lambda_0, \alpha, C_\ast) \epsilon^\frac{\alpha^2}{2}\,,\notag
\end{align}
where we've used again the bound on $|h-H|$ in \eqref{wholebound-2} for $t\in [0,T + T_2]$. Therefore, 
\be \label{boundonH1}
\int_{M_t} H\, d\mu \geq \left(n+\frac{1}{\Lambda_0^2}\right)|M_0| - C_4(n, \Lambda_0, \alpha, C_\ast) 
\epsilon^\frac{\alpha^2}{2} \geq \left(n+\frac{2}{3\Lambda_0^2}\right)|M_0|\,,
\ene
where we've  chosen $\epsilon < \epsilon_{10} = \epsilon_{10}(n, \Lambda_0, \alpha, C_\ast)$ sufficiently small and used the initial condition $\min_{M_0} (H - n) \geq \frac{1}{\Lambda_0^2}$. 

Now applying the bound on $|\nabla H|$ in \eqref{wholebound-2} which holds for all $t\in [0,T + T_2]$, we 
conclude from \eqref{boundonH1} and $|M_t|\leq |M_0|$ that if 
$\epsilon < \epsilon_{11} = \epsilon_{11}(n, \Lambda_0, \alpha, C_\ast)$ is chosen sufficiently small, then for all 
$t\in [0,T + T_2]$, we have
\beq
\min_{M_t} (H - n) \geq \frac{1}{2\Lambda_0^2}\,.
\eeq

Choosing $\epsilon_3 = \min\{\epsilon_4,..., \epsilon_{11}\}>0$, we conclude the proof of the theorem.
\ep

Now we conclude the proof of our main theorem.

\bp (of Theorem \ref{main}) 
In light of Lemma \ref{key}, Lemma \ref{key2} and Theorem \ref{exp1}, by choosing $\Lambda_0$ sufficiently large, we are in position to apply Theorem ~\ref{ext}. 
Thus we can keep extending the VPMCF and estimates (\ref{wholebound-1}) and (\ref{wholebound-2}) for a fixed amount of time depending only on the initial condition. 
Hence the flow \eqref{eq-vp} exists for all time and converges exponentially to a closed umbilic hypersurface in $\H^{n+1}$ by (\ref{wholebound-2}), i.e. an umbilical sphere (\cite{Spi79}).
\ep

\bibliographystyle{amsalpha}
\bibliography{ref-VP}
\end{document}